\documentclass[12pt]{article} 
 \usepackage{mathptmx} 
\usepackage{graphicx,latexsym,amssymb,amsfonts,amsmath,amsxtra,mathdots,units,nicefrac} 
\DeclareMathAlphabet{\mathcal}{OMS}{cmsy}{m}{n}

\newtheorem{theorem}{Theorem} 
\newtheorem{example}{Example} 

\newcommand{\btm}{\begin{theorem}} 
	\newcommand{\etm}{\end{theorem}} 
\newcommand{\beg}{\begin{example}} 
	\newcommand{\eeg}{\end{example}} 
	\newcommand{\bquote}{\vspace{-4pt}\begin{quote}} 
		\newcommand{\equote}{\end{quote}\vspace{-4pt}} 
	\newcommand{\btab}{\vspace{-4pt}\begin{tabbing}}
		\newcommand{\etab}{\end{tabbing}\vspace{-4pt}}
	\newcommand{\beqarray}{\begin{eqnarray*}} 
		\newcommand{\eeqarray}{\end{eqnarray*}} 
\newcommand{\noi}{\noindent}

\newcommand{\q}{\quad}
\newcommand{\h}{\hbox}
 
\newcommand{\R}{\mathbb{R}} 
\newcommand{\N}{\mathbb{N}}
\newcommand{\Z}{\mathbb{Z}}
\newcommand{\Q}{\mathbb{Q}}
\newcommand{\T}{\mathbb{T}}

\newcommand{\sq}{\subseteq} 
\newcommand{\empset}{\varnothing} 

\begin{document}

\centerline{\large\bf Construction and Categoricity of the Real Number System Using Decimals} 

\vspace{.2cm} 
\begin{center} 
Arindama Singh \\ 
Department of Mathematics \\  
Indian Institute of Technology Madras \\ 
Chennai-600036, India \\ 
Email: asingh@iitm.ac.in \\ 
June 27, 2019 
\end{center}

\begin{abstract}
In this expository article, the real numbers are defined as infinite decimals. After defining an ordering relation and the arithmetic operations, it is shown that the set of real numbers is a complete ordered field. It is further shown that any complete ordered field is isomorphic to the constructed set of real numbers. 
\end{abstract}

\noi Keywords: Real numbers, Decimals, Categoricity 

\section{Introduction} 
The modern definition of the system of real numbers posits that it is a complete ordered field. Two questions arise, whether such a theory is consistent, and whether it is categorical. The demand of consistency is met by constructing real numbers from the rational numbers using Dedekind cuts \cite{dede} or by equivalence classes of Cauchy sequences of rationals \cite{cant}. An alternative algebraic construction uses the integers leading to Eudoxus reals; see \cite{str}. Besides these popular constructions of real numbers, a brief survey of other constructions may be found in \cite{weis}.  Once a construction is accepted, the question of categoricity is addressed by proving that a complete ordered field is isomorphic to the system of constructed objects. 

While working with real numbers, one thinks of them as infinite decimal numbers. Indeed,  there have been various attempts at constructing a model for the theory of real numbers by using infinite decimals. The infinite decimals come with the inherent difficulty of defining addition and multiplication. For instance, if two infinite decimals are added, then it is not clear what the $m$th digit of the result will be,  due to the problem of carry digits. In fact, if two infinite decimals are added (or multiplied) and a digit is eventually changed by a carry digit, then it would never change again by another carry digit. Though this fact requires a proof, even accepting it does not lead to determining the $m$th digit of the result. 

An attempt to define real numbers using decimals was taken by K. Weirstrass but it remained somewhat incomplete and unpublished; see \cite{weir}. A later published account may be found in \cite{ritt}, which formalized the notions developed by S. Stevin in \cite{stev}. In this approach real numbers have been defined as formal decimals without bothering whether the rational numbers correspond to a restricted version of these formal decimals; and there was no concern for the modern approach to real numbers via axiomatics. Our approach is similar to this, but we use the decimals as actual decimal numbers and thereby connect to the modern axiomatic definition of real numbers. For this purpose, we assume that we know the following sets along with the usual arithmetic operations and the order $<$: 

\vspace{4pt} 
\noi $\N=\{1,2,3,\ldots\},$ the set of natural numbers. 

\noi $\Z=\{\ldots, -3,-2,-1,0,1,2,3,\ldots\},$ the set of integers. 

\noi $\Q=\{p/q: p\in\Z,~q\in\N,~p,q\h{~have~no~common~factors}\},$ the set of rational numbers. 

\vspace{4pt} 
We know that $\N\subsetneq \Z\subsetneq\Q$ and each rational number can be represented as a recurring decimal number using the {\bf digits} $0,1,2,3,4,5,6,7,8,9$ and the two symbols such as the decimal point $.$ (dot) and the minus sign $-.$ In this representation, any $m\in\Z$ is identified with $\pm m.000\cdots,$ and a decimal number with trailing $9$s is identified with one having trailing $0s.$ that is, $\pm a_0.a_1a_2\cdots a_{k-1}(a_k+1)$ and $\pm a_0.a_1a_2\cdots a_k 9 9 \cdots $ for $0\leq a_k\leq 8,$ are considered equal. Here, $a_0\in\Z$ and $0\leq a_i\leq 9$ for each $i\geq 1,$ and the sign $\pm$ means that either there is no sign or there is a negative sign. The integer $-0$ is again identified with $0.$  A caution about notation: $a_0$ is a finite sequence of digits representing an integer, whereas other $a_i$s are digits. 

With these identifications, the decimal representation of any rational number is unique. In the process, we assume the usual arithmetic of natural numbers, specifically, the algorithms for computing addition and  multiplication of two integers. Further, the notion of less than, $<$, in $\N$  can be specified in terms of digits using the basics such as $0<1 < \cdots < 8<9.$ Consider two unequal natural numbers written in decimal notation, say, $a=a_1a_2\cdots a_k$ and $b=b_1b_2\cdots b_m.$ Then $a<b$ iff one of the following conditions holds: 

\begin{enumerate}  
\item $k<m.$ 
\item $k=m$ and $a_1< b_1.$ 
\item $k=m,~a_1=b_1,\ldots, a_j=b_j$ and $a_{j+1}<b_{j+1}$ for some $j,$ $1\leq j<k.$ 
\end{enumerate}    
Next, the less than relation is extended from $\N$ to $\Z$ by requiring the following conditions: 
\begin{enumerate}  
\item[4.] $0< a_1a_2\cdots a_k.$ 
\item[5.] $-a_1a_2\cdots a_k < b_1b_2\cdots b_m.$
\item[6.] $-a_1a_2\cdots a_k < -b_1b_2\cdots b_m$ iff $b_1b_2\cdots b_m< a_1a_2\cdots a_k.$
\end{enumerate}  

We extend the set of rational numbers to the set of real numbers by formally defining a real number as an infinite decimal, including all terminating, recurring, and non-recurring ones. Our plan is not to extend the rationals to reals but to extend a subset of rational numbers, namely, the terminating decimals to reals. For this purpose, we review certain notions involved with terminating decimals. 

We abbreviate `if and only if' to `iff' as usual. Further, we do not use the pedagogic tool `Definition'; the defined terms are kept bold faced. 

\section{Terminating decimals} 
We assume the usual way of writing integers in decimal notation as explained earlier.  A {\bf terminating decimal} is an expression of the form $a_1a_2\cdots a_k.b_1b_2\cdots b_m$ or $-a_1\cdots a_k.b_1b_2\cdots b_m$ satisfying the following properties: 
\begin{enumerate}  
\item Each $a_i$ and each $b_j$ is a digit. 
\item $a_1\neq 0$ unless $k=1.$ 
\end{enumerate}  
For instance, $51,~51.43,~0.367,~0.0,~0.3600,~-51,~-51.43,~-0.367,~0.0,~-0.3600$ are terminating decimals; and $051.43,~01.3670,~-051.43,~-01.3670$ are not terminating decimals. 

Two terminating decimals are considered {\bf equal} iff either they are identical or one is obtained from the other by adjoining a string of $0$s at the right end. For two terminating decimals $a$ and $b,$ we write $a=b$ when they are equal. That is, given two terminating decimals, we look at how many digits are there in each of them after the decimal point. First, we adjoin necessary number of $0$s at the end of the one which has less number of digits after the decimal point so that both of them now have the same number of digits after the decimal point. Next, we remove the decimal points, and also possible occurrences of the digit $0$ in the beginning, to obtain two integers.  Now, we declare that the two terminating decimals are equal iff the obtained integers are equal. That is, keeping $00\cdots 0a_1a_2\cdots a_k=a_1a_2\cdots a_k$ in the background, we have 
\begin{enumerate}  
\item $a_1a_2\cdots a_k.b_1b_2\cdots b_m =a_1a_2\cdots a_k.b_1b_2\cdots b_m$ 
\item $a_1a_2\cdots a_k.b_1b_2\cdots b_m = a_1a_2\cdots a_k.b_1b_2\cdots b_m 00\cdots 0$ 
\item $a_1a_2\cdots a_k.00\cdots 0 = a_1a_2\cdots a_k. 0=a_1a_2\cdots a_k$ 
\item $-0=0.$ 
\item $-a_1a_2\cdots a_k.b_1b_2\cdots b_m =-c_1c_2\cdots c_k.d_1d_2\cdots d_m$ iff \\ $a_1a_2\cdots a_k.b_1b_2\cdots b_m =c_1c_2\cdots c_k.d_1d_2\cdots d_m.$ 
\end{enumerate}  
This means that formally, we consider $=$ as an equivalence relation on the set of terminating decimals and then update the set of terminating decimals to the equivalence class of $=.$ 

To talk uniformly about terminating decimals, we regard any integer $a\in \Z$ as a terminating decimal $a.b_1b_2\cdots b_m$ with $m=0.$ Notice that the integer $a$ can be in one of the forms $0,$ $a_1a_2\cdots a_k$ or  $-a_1a_2\cdots a_k,$ where $a_1\neq 0.$ This creates no confusion due to the relation of equality as defined above.  We write the set of terminating decimals as $\T.$ This is how $\Z\subsetneq \T.$ 

Two terminating decimals which are not equal, are called unequal, as usual. A nonzero terminating decimal without the minus sign is called a {\bf positive terminating decimal}, and one with the minus sign is called a {\bf negative terminating decimal}.  

Next, two unequal terminating decimals are compared a similar way. First, we adjoin necessary number of $0$s to make the same number of digits after the decimal points in both. Next, we remove the decimal points to obtain two corresponding integers. We then declare that one is {\bf less than} the other iff the corresponding integers are in the same relation. That is, let $a=\pm a_1a_2\cdot a_k.b_1b_2\cdots b_m$ and $c=\pm c_1c_2\cdots c_i.d_1d_2\cdots d_j$ be two unequal terminating decimals. If  $m<j,$ adjoin $j-m$ number of $0$s at the end of $a;$ if $m>j,$ then adjoin $m-j$ number of $0$s at the end of $c$ so that each has equal number of digits after the decimal point. Thus, we take $a'=\pm a_1a_2\cdots a_k.b_1b_2\cdots b_m$ and $c'=\pm c_1c_2\cdots c_i.d_1d_2\cdots d_m.$  We say that $a<c$ iff $a'<c'.$  We read $<$ as `less than', as usual.  

Notice that $<$ is the lexicographic ordering on $\T$, where each negative terminating decimal is less than $0,$ and $0$ is less than each positive terminating decimal. Moreover, $<$ is an extension of the same from $\Z$ to $\T.$  As in $\Z,$ $<$ is a transitive relation on $\T.$  That is, if $a<b$ and $b<c,$ then $a<c.$ Further, we write  $b>a,$ when $a<b.$ As usual, $a\leq b$ abbreviates the phrase `$a<b$ or $a=b$', and $a\geq b$ abbreviates `$a>b$ or $a=b$'. We read these symbols using the phrases, less than, greater than, less than or equal to, or greater than or equal to, as appropriate.   

Next, {\bf addition} and {\bf multiplication} on $\T$ are defined as in usual arithmetic. Due to the equality relation defined earlier (and as we adjoin necessary number of $0$s at the end of one of the terminating decimals) suppose $a=\pm a_1a_2\cdots a_k.b_1b_2\cdots b_m$ and $c=\pm c_1c_2\cdots c_i.d_1d_2\cdots d_m$ be two terminating decimals. Construct integers $a'=\pm a_1a_2\cdots a_kb_1b_2\cdots b_m$ and $c'=\pm c_1c_2\cdots c_id_1d_2\cdots d_m,$ without the decimal points.  Then use the popular addition and multiplication algorithms to obtain $a'+b'$ and $a'\times b'.$ Then $a+b$ is obtained by placing the decimal point preceding the $m$th digit counted from the right of $a'+b'.$ Similarly, $a\times b$ is obtained from $a'\times b'$ by placing the decimal point preceding the $(2m)$th digit in $a'\times b';$ if there do not exist $2m$ digits in $a'\times b',$ then we prefix to it adequate number of $0$s before placing the decimal point. 

It is easy to verify the following properties of addition and multiplication of terminating decimals: 
\begin{enumerate}  
\item $0+a=a$ and $1\times a=a.$ 
\item If $a=c$ and $b=d,$ then $a+c=b+d$ and $a\times c=b\times d.$ 
\item $a+c=c+a,~a\times c=c\times a.$ 
\item $a+(b+c)=(a+b)+c,~a\times (b\times c)=(a\times b)\times c.$ 
\item $a\times (b+c)=(a\times b)+(a\times c).$ 
\end{enumerate}   

Let $A$ and $B$ be nonempty subsets of $\T.$  We define their sum and product as follows: 
$$A+B=\{x+y: x\in A, y\in B\},\q A\times B=\{x\times y: x\in A, y\in B\}.$$ 
Since for terminating decimals $x,y,z,$ we have $x+y=y+x,~x+(y+z)=(x+y)+z,$ $x\times y=y\times x$ and $x\times (y\times z)=(x\times y)\times z,$ we see that for any nonempty subsets $A,B,C,$ of $\T,$   
$$A+B=B+A,~A+(B+C)=(A+B)+C,~A\times B=B\times A,~A\times (B\times C)=(A\times B)\times C.$$ 

\section{Real numbers} 
An {\bf infinite decimal}  is an expression of the form $\pm a_1\cdots a_k.b_1b_2\cdots b_n\cdots,$ where each $a_i$ and $b_j$ is a digit. 

A {\bf real number} is an infinite decimal of the form $\pm a_1\cdots a_k.b_1b_2\cdots b_n\cdots,$ where $a_1\neq 0$ unless $k=1,$ and for no $m,$ all of $b_m, b_{m+1}, b_{m+2}, \ldots $ are $9.$ Again, the symbol $\pm$ says that either the minus sign occurs there or it does not. 

For instance, $21.536800\cdots,~0.536800\cdots,~0.5346464646\cdots,~0.000\cdots,$ \\  $0.101001000100001\cdots,$ $-21.536800\cdots,~-0.536800\cdots,~-0.5346464646\cdots,$ \\ 
$-0.000\cdots$ are  real numbers. In the first and second, the digit $0$ is repeated; in the third, the string $46$ is repeated; in the fourth one, $0$ is repeated, and in the fifth one, nothing is repeated; the others are similar with a minus sign. The expressions $21.53,~021.5399\cdots,~21.5399\cdots,~-21.53,$  $-021.5399\cdots,~-21.5399\cdots$ are not real numbers.  

We say that two real numbers are {\bf equal} iff they are identical. That is, if $a=\pm a_1a_2\cdots a_k.b_1b_2\cdots$ and $c=\pm c_1c_2\cdots c_i.d_1d_2\cdots,$ we say that $a=c$ iff either both start with $-$ or both do not, $k=i,$ and $a_j=c_j$ for each $j\in \{1,\ldots,k\},$ $b_\ell=d_\ell$ for each $\ell\in\N.$ 

Further we identify a real number which ends with a sequence of $0$s with the terminating decimal obtained from it by removing the trailing $0$s. That is, the real number $\pm a_1a_2\cdots a_k.b_1b_2\cdots b_mb_{m+1}\cdots$ with $b_j=0$ for $j>m$ is identified with the terminating decimal  $a_1a_2\cdots a_k.b_1b_2\cdots b_m.$ In this sense of equality, each terminating decimal is regarded as a real number. Thus the real numbers $-0.00\cdots 0\cdots$ and $0.00\cdots 0\cdots$ are each equal to (identified with) the number $0.$ Formally, this identification brings in an equivalence relation and the ensued equivalence classes are the real numbers. 

A nonzero real number is called {\bf positive} if it does not start with the minus sign; and it is called {\bf negative} if it starts with a minus sign. We consider $0$ as  the real number which is neither positive nor negative. A real number is called {\bf non-negative} if it is either $0$ or positive. 

We write the set of all real numbers as $\R.$ Notice that $\T\subsetneq \R.$ 

Two unequal real numbers are compared by defining the relation of $<.$ First, we define $<$ for positive real  numbers. Let $a=a_1a_2\cdots a_k.b_1b_2\cdots$ and $c=c_1c_2\cdots c_i.d_1d_2\cdots$ be unequal positive real numbers. We say that $a<c$ iff one of the following conditions is satisfied: 
\begin{enumerate}  
\item $a_1a_2\cdots a_k<c_1c_2\cdots c_i$ as natural numbers.  
\item $a_1a_2\cdots a_k=c_1c_2\cdots c_i$ and $b_1<d_1.$ 
\item $a_1a_2\cdots a_k=c_1c_2\cdots c_i,~b_1=d_1,\ldots, b_j=d_j$ and $b_{j+1}<d_{j+1}.$ 
\end{enumerate}    
We extend the relation of $<$ to all real numbers in the following manner: 
\begin{enumerate}  
\item If $a$ is negative, then $a<0$ and if $a$ is positive, then $0<a.$ 
\item If $a$ is negative and $c$ is positive, then $a<c.$ 
\item If $a=-a'$ and $c=-c'$ for some positive real numbers $a$ and $c,$ and $c'<a'$ then $a<c.$  
\end{enumerate}  

Again, we say that $b>a$ when $a<b;$ $a\leq b$ abbreviates the phrase `$a<b$ or $a=b$'; and $a\geq b$ abbreviates the phrase `$a>b$ or $a=b$'. We read these symbols in the usual way. It is easy to see that the less than relation is an extension of the same from $\T$ to $\R,$ and the following statements hold for all $a,b,c\in\R$: 
\begin{enumerate}  
\item $a\leq a.$ 
\item If $a\leq b$ and $b\leq a$ then $a=b.$ 
\item If $a\leq b$ and $b\leq c,$ then $a\leq c.$ 
\item $a<b$ or $a=b$ or $b< a.$ 
\item If $a<b$ and $b<c,$ then $a<c.$ 
\item If $a<b,$ then there exists $c\in\T$ such that $a<c<b.$ 
\end{enumerate}  

The last statement above requires some explanation. Suppose $a<b.$ 

First, if $a<0$ and $0<b,$ then take $c=0.0$ 

Second, if $a=0$ and $b$ is positive, then in $b=b_0.b_1b_2\cdots,$ let $m$ be the least index such that $b_m\neq 0.$ Construct $c= 0.b_1b_2\cdots b_{m-1} 01.$ Then $0<c<b.$ 

Third, if $0<a$ and $0<b$, then $a=a_0.a_1a_2\cdots a_ma_{m+1}\cdots$ and $b=b_0.a_1a_2\cdots a_m b_{m+1}\cdots,$ where $a_0$ and $b_0$ are non-negative integers written in decimal notation, and all other $a_i$s are digits. Notice that no infinite decimal ends with repetition of the digit $9.$ Either $a_0<b_0$ or there exists $m\geq 0$ such that $a_i=b_i$ for $0\leq i\leq m$ and $a_{m+1}<b_{m+1}.$ 

If  $a_0<b_0,$ then let $k$ be the least nonzero index such that $a_k<9.$ Then $c=a_0.a_1\cdots a_{k-1}9$ is a required terminating decimal. For instance if $a=1.999988\cdots 8\cdots,$ and $b=2.00\cdots 0\cdots,$ then $c=1.99999;$ and if $a=0.887\cdots 7\cdots$ and $b=5.1\cdots 1\cdots,$ then $c=0.9.$ 

Otherwise, there exists $m\geq 0$ such that $a_i=b_i$ for $0\leq i\leq m$ and $a_{m+1}<b_{m+1}.$   Find out the first occurrence of a digit less than $9$ after the $(m+1)$th digit in $a.$ Change that digit to $9$ and chop-off the rest to obtain $c.$ For instance, if $a=0.12099988\cdots$ and $b=1.12100000\cdots,$ then $m=2$ and the next occurrence of a digit other than $9$ occurs as the $7$th digit after the decimal point; the digit is $8.$ We change this $8$ to $9$ and chop-off the rest to obtain $c=0.1209999.$ Clearly, $a<c<b.$  

Fourth, if both $a$ and $b$ are negative, then let $a'$ and $b'$ be the positive real numbers obtained from $a$ and $b,$ respectively, by removing the minus signs. Now, $b'<a'.$ By the first case, construct $c'$ such that $b'<c'<a'.$ Then $a<-c'<b.$ 

Let $A\sq\R.$  We say that $A$ is {\bf bounded above} iff there exists $b\in\R$ such that $x\leq b$ for each $x\in A.$ In such a case, $b$ is called an {\bf upper bound} for $A.$ For instance, the set of all real numbers less than $2$ is bounded above, for $2$ is one of its upper bounds. Suppose $a_1,a_2,\ldots,a_k$ are fixed digits. Then the set of all real numbers of the form $\pm a_1a_2\cdots a_k.b_1b_2\cdots$ is bounded above since $a_1a_2\cdots a_k0$ is one of its upper bounds. The set $\N=\{1,2,3,\ldots\}$ is not bounded above.

Clearly, if $b$ is an upper bound of a set $A$ of real numbers, then any real number $c>b$ is also an upper bound of $A.$ However, a real number smaller than $b$ need not be an upper bound of $A.$ For instance, The infinite set $\{0.9,~0.99,~0.999,~\ldots\}$ of terminating decimals has an upper bound $1.$ But no real number smaller than $1$ is an upper bound of this set. Reason? Any real number less than $1=1.00\cdots$ is in the form $0.b_1b_2\cdots.$ If it is an upper bound of the set, no $b_i$ is less than $9.$ Then the infinite decimal must be $0.999\cdots.$ However, it is not a real number! 

If $a=a_1a_2\cdots a_k.b_1b_2\cdots$ is a non-negative real number, write $[a]=a_1a_2\cdots a_k;$ and   \\ 
\indent if $a=-a_1a_2\cdots a_k.b_1b_2\cdots$ is a negative real number, write $[a]=-a_1a_2\cdots a_k-1.$ \\ 
This is how for each real number $a,$ there exists an integer $[a]$ such that $[a]\leq a < [a]+1.$ Here,  we use the known operation of addition of integers by writing an integer as a concatenation of digits with or without the minus sign. Such an integer $[a],$ which is equal to the real number $[a].00\cdots,$ is called  the {\bf integral part} of $a.$ 

For a nonempty subset $A$ of $\R$ which is bounded above, we define its {\bf supremum}, denoted by $\sup(A),$  by considering two cases. \\[4pt] 
{\em Case 1}: Suppose $A$ has at least one non-negative real number. Let $A_0$ be the set of integral parts of real numbers in $A.$ Notice that $A_0$ contains at least one non-negative integer. Since $A$ is bounded above, there exists a maximum of $A_0.$ Write the maximum as $m.$ Now, $m$ is a non-negative integer.  In $A,$ there are real numbers whose integral part is $m.$  Let $A_1$ be the set of all real numbers in $A$ whose integral part is $m.$  Each real number in $A_1$ is in the form $m.b_1b_2\cdots.$ Look at the first decimal digit $b_1$ in all such real numbers. There exists one maximum digit, say $m_1.$ Let $A_2$ be the set of all real numbers in $A_1$ which are in the form $m.m_1c_2c_3\cdots.$ Then choose the maximum possible digit $c_2$ among these real numbers, and call it $m_2.$ 

Continuing in this fashion, we obtain a sequence of digits $m_1,m_2,m_3,\ldots.$ Then we form the infinite decimal $m.m_1m_2m_3\cdots.$ This infinite decimal may or may not be a real number. We consider two cases. \\[4pt] 
{\em Case }(1A): Suppose there is no such $i$ such that $m_j=9$ for all $j\geq i.$ That is, the infinite decimal does not end with a sequence of $9$s. Then, $m.m_1m_2m_3\cdots$ is a real number; and we call this real number as $\sup(A).$ \\[4pt] 
{\em Case} (1B): Suppose there exists $i$ such that $m_j=9$ for all $j\geq i.$ That is, the $i$th digit onwards all digits are equal to $9.$ Then consider the terminating decimal $m.m_1m_2m_3\cdots m_i$ and add to it $0.00\cdots 1$ having the $i$th digit after the decimal as $1$ and all previous digits $0.$ The result is a terminating decimal. This terminating decimal is a real number, which we call $\sup(A).$ \\[4pt] 
{\em Case 2}: Suppose that all elements in $A$ are negative. Let $A'=\{x: -x\in A\}.$ That is, we consider real numbers in $A$ after removing the minus signs. We follow the construction in Case 1, but this time by taking minimum instead of maximum everywhere, and then put a minus sign in the result. The details follow. 

Notice that if $-u$ is an an upper bound of $A,$ then $u$ is less than or equal to each element of $A'.$ Let $A'_0$ be the set of integral parts of (non-negative) real numbers in $A'.$  There exists a smallest non-negative integer among them, that is, the minimum of $A'_0.$ Write the minimum as $n.$ In $A',$ there are real numbers whose integral part is $n.$  Let $A'_1$ be the set of all real numbers in $A'$ whose integral part is $n.$  Each real number in $A'_1$ is in the form $n.d_1d_2\cdots.$ Look at the first decimal digit $d_1$ in all such real numbers. There exists one minimum digit, say $n_1.$ Let $A'_2$ be the set of all real numbers in $A'_1$ which are in the form $n.n_1e_2e_3\cdots.$ Then choose the minimum possible digit $e_2$ among these real numbers, and call it $n_2.$ 

Continuing in this fashion, we obtain a sequence of digits $n_1,n_2,n_3,\ldots.$ Then we form the infinite decimal $n.n_1n_2n_3\cdots.$ This infinite decimal may or may not be a real number. We consider two cases. \\[4pt] 
{\em Case }(2A): Suppose there is no such $i$ such that $n_j=9$ for all $j\geq i.$ That is, the infinite decimal does not end with a sequence of $9$s. Then, $n.n_1n_2n_3\cdots$ is a real number. We call the real number $-n.n_1n_2n_3\cdots$ as $\sup(A).$ \\[4pt] 
{\em Case} (2B): Suppose there exists $i$ such that $n_j=9$ for all $j\geq i.$ That is, $i$th digit onwards all digits are equal to $9.$ Then consider the terminating decimal $n.n_1n_2n_3\cdots n_i$ and add to it $0.00\cdots 1$ having the $i$th digit after the decimal as $1$ and all previous digits $0.$ The result is a terminating decimal. Prefix this terminating decimal with a minus sign, and cal, it $\sup(A).$    

\beg 
Compute the supremum of the following subsets of $\:\R$: \\ 
$A=\{1.0, ~2.12,~1.11\cdots 1\cdots, ~2.12011\cdots 1\cdots,~1.1201011\cdots 1\cdots\},$  \\  $B=\{0.9,~0.99,~0.19,~0.991,~0.9991,~0.99991, ~\cdots, ~0.99\cdots 9 1,~\cdots\},$ \\ 
$C=\{-1.0,~-0.9,~-0.99,~-0.19,~-0.991,~-0.9991,~-0.99991, \cdots, ~-0.99\cdots 9 1,~\cdots\}.$  \\ 
$D=\{-0.1,~-0.01,~-0.001,~-0.0001,~\cdots,~-0.00\cdots 01,~\cdots\}.$  
\eeg 
\noi For $A,$ the set of integral parts is $A_1=\{1,2\}.$ So $m_1=2.$ Next, $A_2=\{2.12,2.12011\cdots 1\cdots\}.$ So, $m_2=1.$ $A_3=\{2.12,2.12011\cdots 1\cdots\}.$ Then $m_3=2.$ Next, $A_4=\{2.12,2.12011\cdots 1\cdots\}.$ We use the equality of real numbers and see that $2.12=2.1200\cdots 0\cdots.$ So, $m_4=0.$  Then $A_5=\{2.12,2.12011\cdots 1\cdots\},$ and $m_5=1.$ Next, $A_6=\{2.12011\cdots 1\cdots\}.$ So, $m_6=1.$ This point onwards we get each $m_j$ to be equal to $1.$ Hence  $\sup(A)=2.12011\cdots 1\cdots.$ \\[4pt] 
For $B$, $m=0,~m_1=9,~m_2=9,$ and all succeeding digits are $9.$ Then we have the infinite decimal $0.999\cdots 9\cdots.$ We see that $i=1.$ So, we consider the decimal $0.9$ and add to it $0.1$ to obtain $\sup(B)=1.$ \\[4pt] 
For $C,$ the set $C'=\{1.0,~0.9,~0.99,~0.19,~0.991,~0.9991,~0.99991, \cdots, ~0.99\cdots 9 1,~\cdots\}.$ Then $n=0,~C'_1=\{0.9,~0.99,~0.19,~0.991,~0.9991,~0.99991, \cdots, ~0.99\cdots 9 1,~\cdots\}.$ So that $n_1=1.$ Now, $C'_2=\{0.19\},$ so that $n_2=9,~n_3=0,\ldots.$ So, $\sup(C)=-0.19.$ \\[4pt] 
For $D,$ the set $D'=\{0.1,~0.01,~0.001,~\cdots,~0.00\cdots 01,~\cdots\}.$ Then $n_1=0,$ $n_2=0,~\ldots.$ Therefore, $\sup(D)=-0.00\cdots 0\cdots =0.$ \hfill $\diamond$ 

As we see, each set of real numbers which is bounded above has a unique supremum, which is also a real number.  We give a characterization of the supremum. 

\btm\label{supch}  
Let $A$ be a nonempty subset of $\:\R$ which is bounded above and let $s\in \R.$ Then, $s=\sup(A)$ iff $s$ is an upper bound of $A$ and one of the following conditions is satisfied: 
\begin{enumerate}  
\item If $t\in \R$ is an upper bound of $A,$ then $s\leq t.$ 
\item If $p\in\R$ and $p<s,$ then there exists $q\in A$ such that $p<q\leq s.$ 
\end{enumerate}  
\etm 

\noi{\em Proof.} Suppose $S$ is a finite set of negative integers, say, $S=\{-n_1,-n_2,\ldots,-n_k\}.$ Let $S'=\{n_1,n_2,\ldots,n_k\}.$ Clearly, $-n_1$ is the minimum of $S$ iff $n_1$ is the maximum of $S'.$ Moreover, suppose a set $X$ of real numbers contains at least one non-negative real number. Let $Y$ be the set of all non-negative real numbers that belong to $X.$ Then clearly, $\sup(X)=\sup(Y).$ Therefore, all the cases discussed in the construction of the supremum reduce to the case where a set has only non-negative real numbers. So, without loss of generality, assume that $A$ is a set of non-negative real numbers. 

From our construction it follows that  $\sup(A)$ is an upper bound of $A.$ \\[4pt] 
(1) Let $t$ be an upper bound of $A$ but $t\neq \sup(A).$ First, the integral part of each element of $A$ is less than or equal to that of $t.$ Next, for each $x\in A,$ the $m$th digit of $x$ is less than or equal to the $m$th digit of $t.$  Since $t\neq\sup(A),$ some $m$th digit of $t$ is greater than the $m$th digits of all real numbers in $A.$ However, the $m$th digit of $\sup(A)$ is the $m$th digit of some real number in $A.$ Therefore, the $m$th digit of $\sup(A)$ is less than the $m$th digit of $t.$ So, $\sup(A)<t.$ Therefore, $\sup(A)\leq t.$ 

Conversely, suppose that $s$ is an upper bound of $A$ and that $s$ is less than or equal to every upper bound of $A.$ As $\sup(A)$ is an upper bound of $A,$ $s\leq \sup(A).$ However, by what we have proved in the above paragraph, $\sup(A)\leq s.$ So, $s=\sup(A).$ \\[4pt] 
(2) Let $p$ be a real number such that $p<\sup(A).$ If no real number in $A$ is larger than $p,$ then $p$ is an upper bound of $A.$ By (1), $\sup(A)\leq p$ which is a contradiction. 

Conversely, suppose that $s$ is an upper bound of $A$ and that if $p$ is a real number with $p<s,$ then there exists a real number $q\in A$ such that $p<q.$ If $s\neq\sup(A),$ we have $\sup(A)<s.$ Then with $p=\sup(A),$ we have a real number $q\in A$ such that $\sup(A)<q.$ This is a contradiction to the fact that $\sup(A)$ is an upper bound of $A.$ Once $p<q,~q\in A$ and $s$ is an upper bound of $A,$ it follows that $p<q\leq s.$  \hfill $\Box$

The above theorem says that if $A$ is a nonempty set of real numbers bounded above, then $\sup(A)$ is the least of all upper bounds of $A.$ Some useful properties of the supremum are proved in the following theorem. 

\btm\label{forall-exists-sup}  
Let $A$ be a nonempty subset of $\:\R$ which is bounded above. 
\begin{enumerate}  
\item Let $B$ be a nonempty subset of $\R$ that is bounded above. If for each $x\in A,$ there exists $y\in B$ such that $x\leq y,$ then $\sup(A)\leq \sup(B).$ In particular, if $A\sq B,$ then $\sup(A)\leq \sup(B).$ 
\item If $a,b\in\R,~a<b,$ then $\sup\{x\in\R: x<b\}=b=\sup\{x\in\R: a<x<b\}.$ 
\item If $c\in\R,$ then $c=\sup\{x\in\T: x<c\}.$  
\item $\sup\{\sup\{t\in\T: t<x\}: x\in A\}=\sup(A)=\sup\{t\in\T: t<\sup(A)\}.$ 
\end{enumerate}   
\etm 

\noi{\em Proof.} (1) Due to our assumption, for each $x\in A,~x\leq \sup(B).$ Hence $\sup(B)$ is an upper bound of $A.$ Since $\sup(A)$ is less than or equal to each upper bound of $A,$ we have $\sup(A)\leq \sup(B).$ \\[4pt] 
(2) Let $B=\{x\in\R: x<b\}.$ Clearly, $b$ is an upper bound of $B.$ If $u\in \R$ is such that $u<b,$ then there exists $v\in\T$ such that $u<v<b.$ Now, $v\in B.$ Hence $b=\sup(B).$ 

For the second equality, let $a\in\R$ and $a<b.$ Write $C=\{x\in\R: a<x<b\}.$  Notice that for each $x\in C,$ there exists $y\in B$ such that $x\leq y;$ and for each $y\in B,$ there exists $x\in C$ such that $y\leq x.$ Hence by (1), $\sup(C)=\sup(B)=b.$ \\[4pt] 
(3) Let $c\in\R,~D=\{x\in\R: x<c\}$ and let $E=\{x\in \T: x<c\}.$ If $x\in D,$ then there exists $y\in\T$ such that $x<y<c.$ That is, $x<y$ for some $y\in E.$ Moreover, $E\sq D.$ Hence by (1)-(2),  $\sup(E)=\sup(D)=c.$ \\[4pt] 
(4) By (3) and (2), $\sup\{\sup\{t\in\T: t<x\}: x\in A\}=\sup\{x: x\in A\}=\sup(A)=\sup\{t\in\T: t<\sup(A)\}.$  
\hfill $\Box$ 

\pagebreak %\vspace{6pt} 
We define {\bf addition} and {\bf multiplication} of two real numbers $a$ and $b$ as follows. 
\begin{enumerate}  
\item If $a=0,$ then $a+b=b$ and $a\times b=0.$ 
\item If $b=0,$ then $a+b=a$ and $a\times b=0.$ 
\item If $a\neq 0$ and $b\neq 0,$ then $a+b=\sup(A+B),$ where $A=\{x\in\T: x < a\}$ and $B=\{x\in\T:x < b\}.$ 
\item If $0<a$ and $0<b,$ then $a\times b = \sup(A\times B),$ where $A=\{x\in\T: 0<x < a\}$ and $B=\{x\in\T:0<x < b\}.$ 
\item If $0<a$ and $b<0,$ then $b=-y,~0<y$ so that $a\times b = -(a\times y).$ 
\item If $a<0$ and $0<b,$ then $a=-x,~0<x$ so that $a\times b=-(x\times b).$ 
\item If $a<0$ and $b<0,$ then $a=-x,~b=-y,0<x,~0<y$ so that $a\times b=x\times y.$  
\end{enumerate}  

Notice that in the cases (3)-(4), the sets $\{x+y: x\in A, y\in B\}$ and $\{x\times y: x\in A, y\in B\}$  are bounded above. For instance, if $0<a=a_1a_2\cdots a_k.c_1c_2\cdots$ and $0<b=b_1b_2\cdots b_i.d_1d_2\cdots,$ then upper bounds of $A+B$ and $A\times B$ can be given as $a_1a_2\cdots a_k+b_1b_2\cdots b_i+2$ and $(a_1a_2\cdots a_k+1)\times (b_1b_2\cdots b_i +1),$ respectively. Similar constructions work for other possibilities of $a$ and $b.$  Hence addition and multiplication of real numbers are well defined. 

\btm\label{plusdot} 
Let $A$ be a nonempty subset of $\:\R$ that is bounded above and let $b\in\R.$ Then 
$~\sup\{b+x: x\in A\}=b+\sup(A)~\h{~and~}~\sup\{b\times x: x\in A\}=b\times \sup(A).$ 
\etm 

\noi {\em Proof.} When $b=0,$ the result is obvious. For $0<b,$ using the definition of $+,~\times$ and Theorem~\ref{forall-exists-sup}(4), we obtain 
\beqarray 
\sup\{b+x: x\in A\} & = & \sup\{\sup\{y+t: y\in\T,~t\in\T,~y<b,~t<x\}: x\in A\} \\ 
& = & \sup\{y+t: y\in\T,~t\in\T,~y<b,~t<\sup(A)\} = b+ \sup(A). \\ 
\sup\{b\times x: x\in A\} & = & \sup\{\sup\{y\times t: y\in\T,~t\in\T,~y<b,~t<x\}: x\in A\} \\ 
& = & \sup\{y\times t: y\in\T,~t\in\T,~y<b,~t<\sup(A)\} = b\times \sup(A). 
\eeqarray 
Similarly, the case $b<0$ is tackled. \hfill $\Box$

In particular, if $A\sq\T$ that is bounded above, then the conclusion of the above theorem also holds. 

We now proceed to prove the required properties of real numbers. 

\btm\label{axioms} Let  $x,y,z\in \R$ and let $S$ be a nonempty subset of $\R.$ Then the following are true:  

\begin{enumerate}   
\item $x+y=y+x.$ 
\item $(x+y)+z= x+(y+z).$  
\item $0+x=x.$ 
\item Corresponding to $x$ there exists $u\in\R$ such that $x+u=0.$ 
\item $x\times y=y\times x.$ 
\item $(x\times y)\times z = x\times (y\times z).$ 
\item $1\times x = x.$ 
\item Corresponding to $x\neq 0,$ there exists $w\in \R$ such that $x\times w=1.$ 
\item $x\times (y+z)=(x\times y)+(x\times z).$ 
\item Exactly one of the conditions $x<y$ or $x=y$ or $y<x$ is true. 
\item If $x<y$ and $y<z,$ then $x<z.$ 
\item If $x<y,$ then $x+z<y+z.$ 
\item If $x<y$ and $0<z,$ then $z\times x< z\times y.$ 
\item If there exists $b\in S$ with $x\leq b$ for each $x\in S$, then there exists $s\in \R$ such that for each $x\in S,~x\leq s$; and for each $u\in\R,$ if $x\leq u,$ then $s\leq u.$ 
\end{enumerate}        
\etm 

\noi {\em Proof.} Let $A,B,C$ be the corresponding sets of terminating decimals less than $x,y,z,$ respectively.  By definition, $x+y=\sup(A+B)$ and etc. \\[4pt] 
(1) Since $A+B=B+A,$  we have $x+y=\sup(A+B)=\sup(B+A)=y+x.$ \\[4pt] 
(2) On the contrary, suppose $x+(y+z)<(x+y)+z.$ Let $d\in\T$ be such that $x+(y+z)<d<(x+y)+z.$ Now, $d<(x+y)+z$ implies that there exist $p,q\in\T$ such that $d<p+q,~p< x+y,~q < z.$ \\ 
Again, $p<x+y$ implies that  there exist $r,s\in\T$ such that $p<r+s,~r<x,~s<y.$ \\ 
Then $d<(r+s)+q=r+(s+q).$ But $s\in\T,~s<y,$ and $q\in\T$ with $q<z.$ So, $s+q <y+z.$ Since $r\in\T$ and $r<x,$ we see that $r+(s+q)<x+(y+z).$ Hence $d<x+(y+z).$ This is a contradiction since $x+(y+z)<d.$ \\ 
\indent Similarly, $(x+y)+z<x+(y+z)$ leads to a contradiction. Therefore, $x+(y+z)=(x+y)+z.$ \\[4pt] 
(3) It follows by definition. \\[4pt] 
(4) We show that $x+(-x)=0.$ Let $A=\{s\in\T: s<x\}$ and let $B=\{t\in\T: t<-x\}.$ Clearly, $A,~B$ and $A+B$ are nonempty, and $A,~B$ are bounded above. Now, if $s\in A$ and $t\in B,$ then $s<x<-t$ so that $s+t<t+(-t)=0.$ That is, $A+B$ is bounded above by $0.$ Also, if $r\in\T$ and $r<0,$ then with any $s\in A$ we have $r+(-s)< 0+(-s)=-s <-x.$ That is, with $t=r+(-s),$ we see that $s+t=r.$ It implies that each $r\in\T$ with $r<0$ belongs to $A+B.$ Hence, $A+B=\{y\in\T: y<0\}.$ By Theorem~\ref{forall-exists-sup}(4), $x+(-x)=\sup(A+B)=0.$ \\[4pt] 
(5)-(6) Proofs of these are similar to (1)-(2). \\[4pt] 
(7) By Theorem~\ref{plusdot}, $1\times x=\sup\{1\times t: t\in\T,~t<x\}=\sup\{t:t\in\T,~t<x\}=x.$ \\[4pt] 
(8) Let $0<x.$ Write $A=\{t\in\T: x\times t<1\}.$ If $x=a_1a_2\cdots a_k.b_1b_2\cdots,$ then take $u=0.00\cdots 01,$ with $k+1$ number of $0$s. We see that $x\times u<1.$ So, $A\neq\empset.$ Further, if $x\geq 1,$  then $1$ is an upper bound of $A.$ If $0<x<1,$ then suppose there are $m$ number of $0$s in $x$ that immediately follow the decimal point. Then $10^{m+1},$ that is, the natural number $1$ followed by $m+1$ number of $0$s is an upper bound of $A.$ Hence $A$ is bounded above. Thus, it has a supremum. Let $w=\sup(A).$ Notice that $0<w.$ By Theorems~\ref{plusdot} and \ref{forall-exists-sup} (2), $x\times w = x\times \sup(A)= \sup\{x\times t: t\in A\} = \sup\{x\times t: x\times t<1\}=1.$  \\ 
\indent Next, if $x<0,$ let $x=-y,$ where $0<y.$ By what we have just proved, there exists $z\in\R,~0<z$ such that $y\times z=1.$ Then $x\times (-z)=(-y)\times (-z)=y\times z=1.$   \\[4pt] 
(9) This follows from a similar equality for terminating decimals by taking supremum. \\[4pt] 
(10)-(11) These properties follow directly from the definition of $<.$ \\[4pt] 
(12) Let $x<y.$ Write $A=\{r\in \T: r<x\},~B=\{s\in T: s<y\}$ and $C=\{t\in T: t<z\}.$ Now $A\sq B$ so that $A+C\sq B+C.$ By Theorem~\ref{forall-exists-sup} (1), $x+z=\sup(A+C)\leq \sup(B+C)=y+z.$ If $x+z=y+z,$ then by (1)-(4), we get $x=y,$ which is not possible. Therefore, $x+z<y+z.$ \\[4pt] 
(13)  Let $x<y$ and $0<z.$ Write $A=\{r\in \T: r<x\},~B=\{s\in T: s<y\}$ and $C=\{t\in\T: 0<t<z\}.$ Let $u\in A\times C.$ Then $u=r\times t,$ where $r,t\in\T,$ $r<x<y$ and $0<t<z.$ Thus $u\in B\times C.$ Hence $x\times z =\sup(A\times C)\leq \sup(B\times C) = y\times z.$ If $x\times z= y\times z,$ then using (6)-(8) we would get $x=y,$ which is a contradiction. Therefore, $x\times z=y\times z.$ \\[4pt] 
(14) Due to Theorem~\ref{supch}, this statement may be rephrased as:  if a nonempty subset $S$ of $\R$ is bounded above, then $s=\sup(S)\in\R.$ This is obvious from the construction of supremum. \hfill $\Box$

\section{Axiomatic definition of real numbers} 
The modern definition of $\R$ uses the properties of $\R$ as axioms. In such a treatment, we define ${\cal R}$ as a set containing at least two distinct symbols $0$ and $1,$ where two binary operations, denoted as $+$ and $\times,$ and a binary relation called less than, denoted by $<,$ are assumed given; and $+,~\times, <$ are such that the fourteen conditions specified in Theorem~\ref{axioms} with $\R$ replaced by ${\cal R}$ are satisfied. These conditions are called axioms of ${\cal R}.$ 

It can be shown that an element $u\in {\cal R}$ corresponding to any $x\in {\cal R}$ in Axiom~(4) is unique, and we write it as $-x.$ Similarly, it can be shown that an element $w\in{\cal R}$ corresponding to a nonzero $x\in {\cal R}$ in Axiom~(8) is unique, and we write it as $1/x.$ 

In fact, any set with at least two elements, called $0$ and $1,$ satisfying the first nine axioms is called a {\bf field}.  A field where a binary relation, denoted as $<$ is defined, satisfying the axioms (10)-(13) is called an {\bf ordered field}. In any ordered field, an element $x$  is called {\bf positive} if $0<x;$ and $x$ is called {\bf negative} if $x<0.$ The element $0$ is neither positive nor negative. Further, it can be shown that $0<1.$ 

Notice that both $\R$ and $\Q$ are ordered fields. To distinguish $\R$ from $\Q$ the fourteenth condition is assumed. The fourteenth axiom of ${\cal R}$ concerns not every element of ${\cal R}$ but every subset of ${\cal R}.$ This condition is worth reformulating introducing new simpler phrases. Suppose $S$ is a nonempty subset of ${\cal R}.$ A number $b\in {\cal R}$ is called an {\bf upper bound of} $S$ iff every element of $S$ is less than or equal to $b.$ A number $s\in {\cal R}$ is called a {\bf least upper bound of} $S$ in ${\cal R}$ iff $s$ is an upper bound of $S$ and no upper bound of $S$ is less than $s.$ We abbreviate the phrase `least upper bound' to {\bf lub}.  Axiom~(14) asserts that every subset of ${\cal R}$ having an upper bound has an lub in ${\cal R}.$ It is called the completeness axiom. 

The following characterization of the ${\rm lub}$ holds in ${\cal R}.$ 

\btm 
Let $S$ be a nonempty subset of $\:{\cal R}$ which is bounded above and let $s\in R.$ Then, $s={\rm lub}(S)$ iff $s$ is an upper bound of $S$, and if $p\in {\cal R}$ and $p<s,$ then there exists $q\in S$ such that $p<q\leq s.$ 
\etm 

\noi {\em Proof.} Let $s={\rm lub}(S).$ Then $s$ is an upper bound of $S.$ Let $p\in {\cal R}$ and $p<s.$ Then $p$ is not an upper bound of $S.$ So, there exists $q\in S$ such that $p<q.$ As $s$ is an upper bound of $S,$ we also have $q\leq s.$ That is, such a $q\in S$ satisfies $p<q\leq s.$ 

Conversely, suppose that $s$ is an upper bound of $S$ and that if $p\in {\cal R}$ satisfies $p<s,$ then there exists $q\in S$ such that $p<q\leq s.$ Now, if $r\in {\cal R}$ satisfies $r<s,$ then $r$ is not an upper bound of $S.$ Hence $s={\rm lub}(S).$ \hfill $\Box$

It follows that lub of a set, which is bounded above, is unique. The axioms assert that ${\cal R}$ is a complete ordered field. We give an application of the completeness axiom of ${\cal R},$ which will be useful later.  In what follows, we will use the commonly accepted abbreviations such as writing $x\times y$ as $xy,$ $x\times x$ as $x^2,$ $x+(-y)$ as $x-y,$ $x\times (1/y)$ as $x/y,$ and the precedence rules, where $(x\times y)+(z\times w)$ is abbreviated to $xy+zw,$ and etc. 

\btm\label{sqroot}  
Let $r\in {\cal R}$ be positive. Then there exists a unique positive $s\in {\cal R}$ such that $s^2=r.$ 
\etm 

\noi {\em Proof.} If $r=1,$ then take $s=1$ so that $s^2=1.$ If $r>1,$ then $1/r<1.$ If there exists $t>0$ such that $t^2=1/r,$ then $s=1/t$ satisfies $s^2=r.$ So, without loss of generality, suppose $0<r<1.$ 

Write $A=\{x\in {\cal R}: 0<x^2\leq r\}.$ Now, $r\in A$ and $A$ is bounded above by $1.$ So, let $s={\rm lub}(A).$ Notice that $s\leq r<1.$ 

If $r<s^2,$ then take $\epsilon=s^2/2.$ Now, $\epsilon>0,$ and $r=s^2-2\epsilon < s^2-2s\epsilon < (s-\epsilon) ^2.$ That is, $s-\epsilon$ is an upper bound of $A.$ This contradicts the fact that $s={\rm lub}(A).$ 

If $s^2<r,$ then take $\epsilon=\min\Big\{\dfrac{r-s^2}{2s+1},1\Big\}.$ Now, $s+\epsilon>s$ and 
$$(s+\epsilon) ^2=s^2+2s\epsilon+\epsilon^2 \leq s^2+2s\epsilon+\epsilon = s^2+\epsilon(2s+1)\leq s^2+ \dfrac{r-s^2}{2s+1} (2s+1)= r.$$ 
That is, $s+\epsilon \in A.$ It contradicts the fact that $s={\rm lub}(A).$  

Therefore, $s^2=r.$ 

To show uniqueness of such a positive $s,$ suppose $t\in {\cal R},$ $0<t$ and $t^2=r.$ If $s<t,$ then $0<s^2<s t < t^2$ contradicts the assumption that $s^2=r=t^2.$ Similarly, $t<s$ is not possible. That is, $s=t.$ \hfill $\Box$

We write the unique positive $s\in {\cal R}$ that satisfies $s^2=r$ as $\sqrt{r},$ and call it the {\bf positive square root} of the given positive element $r\in {\cal R}.$ 

We have not assumed that our known sets of numbers such as $\N,$ $\Z,$ or $\Q$ are contained in ${\cal R}.$ Reason? They are already inside ${\cal R}$ in some sense. How? 

Let ${\cal R}$ be a complete ordered field. Let us look at $\N$ first. The number $1$ of $\N$ is in ${\cal R}.$ Of course, the $1$ inside ${\cal R}$ is just a symbol, like $1$ inside $\N.$ Well, we make a correspondence of $1$ inside $\N$ to $1$ inside ${\cal R}.$ Now, $1+1$ inside $\N$ corresponds to $1+1$ inside ${\cal R}.$ Given a sum of $n$ number of $1$s inside $\N$ now corresponds to that inside ${\cal R}.$ In this sense, $\N$ is inside ${\cal R},$ that is, a copy of $\N$ is inside $R.$ But is the induction principle of $\N$ available to this copy of $\N$ inside ${\cal R}$? The answer is affirmative, and we show it as follows. 

Let $S\sq R.$ We call $S$ an {\em inductive subset} of ${\cal R}$ iff $1\in S$ and for each $x\in S,$ $x+1\in S.$ For example, ${\cal R}$ is an inductive subset of itself. Now, let ${\cal N}$ be the intersection of all inductive subsets of ${\cal R}.$ It follows that if $S$ is any inductive subset of ${\cal R}$, then ${\cal N}$ is a subset of $S.$ Further, let $A\sq {\cal N}.$ If $A$ is an inductive set in ${\cal R},$ then using what we have just proved, we obtain $A={\cal N}.$ Thus, we have proved the following: 
\bquote 
Let $S$ be a subset of ${\cal N}$ satisfying (i) $1\in S;$ and (ii) if $x\in S,$ then $x+1\in S.$ Then $S={\cal N}.$ 
\equote 
This is exactly the induction principle in $\N.$ Therefore, a proper copy of $\N$ inside ${\cal R}$ is this set ${\cal N},$ the intersection of all inductive subsets of ${\cal R}.$ 

Recall that the induction principle implies the well ordering principle, which states that 
\bquote 
Each nonempty subset of ${\cal N}$ has a least element. 
\equote 

Then ${\cal N}$ is extended to ${\cal Z}$ by taking ${\cal Z}={\cal N}\cup \{0\}\cup \{-n: n\in {\cal N}\}$ with the relation of $<$, the operations $+$ and $\times $ as in ${\cal R}.$ Next, ${\cal Q}$ is obtained from ${\cal Z}$ by taking ${\cal Q}=\{p/q: p\in Z, q\in N,\h{~and~}p,q\h{~have~no~common~factors}\}.$  Now, the sets ${\cal N},~{\cal Z},~{\cal Q}$ are copies of $\N,~\Z,~\Q,$ inside ${\cal R}$, respectively. We show that ${\cal Q}\neq {\cal R}$ by using the axioms of ${\cal R}$ and the construction of ${\cal Q}$ inside ${\cal R}.$ 

\btm\label{qincomp}   
${\cal Q}$ is an ordered field, but it is not complete. 
\etm  

\noi {\em Proof.} It is obvious that ${\cal Q}$ is an ordered field. On the contrary, suppose that ${\cal Q}$ is complete. Notice that Theorem~\ref{sqroot} holds in every complete ordered field, and in particular, for ${\cal Q}.$ Now that $2\in {\cal Q},$ we see that there exists $s\in {\cal Q}$ such that $s^2=2.$ Then $s=\frac{p}{q},$ where $p\in {\cal Z}$ and $q\in {\cal N}$ do not have any common factor. Then $s^2=2$ implies that $p^2=2q^2.$ Since $p^2$ is an even number and a square, it follows that $p$ is a multiple of $4$ (by induction). In that case, $p=2m$ for some $m\in {\cal Z}.$ That is, $4m^2=2q^2.$ It gives $2m^2=q^2.$ With a similar argument, we see that $q=2n$ for some $n\in {\cal N}.$ This is a contradiction since $p$ and $q$ do not have any common factor.  \hfill $\Box$

The above shows, essentially, that $\sqrt{2}\not\in {\cal Q}.$ 

Notice that $\N$ and its copy ${\cal N}$ inside ${\cal R}$ are isomorphic, in the sense that if $x,y\in\N$ are written respectively as $x',y'\in {\cal N},$ then $x+y$ and $x\times y$ are renamed as $x'+y'$ and $x'\times y',$ respectively. Here, the $+$ and $\times $ are the same operations as they are given in ${\cal R}.$ Same way, $\Z$ and its copy ${\cal Z}$ inside ${\cal R}$ are isomorphic. It can be shown that $\Q$ is isomorphic to ${\cal Q}$ inside ${\cal R}.$ Further, the isomorphism of $\Q$ with ${\cal Q}$ is an extension of the said isomorphisms of $\N$ with ${\cal N}$ and of $\Z$ with ${\cal Z}.$ In addition, this isomorphism of $\Q$ to ${\cal Q}$ also preserves the relation of $<.$ Reason is, the same conditions (10)-(13) hold in $\Q$ as well as in ${\cal Q}.$ 

\section{Uniqueness of ${\cal R}$} 
Notice that $\R$ satisfies all the axioms of ${\cal R}.$ Hence $\R$ is a model of the theory ${\cal R}.$ Our construction of $\R$ assumes the constructions of $\N$ and $\T.$ Again, $\T$ is constructed from $\N.$ Hence it shows that if $\N$ is a consistent theory, then so is ${\cal R}.$ We will show that ${\cal R}$ is a categorical theory, in the sense that the axioms (1)-(14) define the object ${\cal R}$ uniquely up to an isomorphism, by proving that ${\cal R}$ is isomorphic to $\R.$ In fact, the decimal representation of any element of ${\cal R}$ provides such an isomorphism. For a complete treatment, we give the details in the following. 

\btm 
(Archimedean Principle): Let $x,y\in {\cal R}$ with $x>0.$ Then there exists $n\in {\cal N}$ such that $nx>y.$ 
\etm 

\noi {\em Proof}. Suppose $nx\leq y$ for all $n\in {\cal N}.$ Then $y$ is an upper bound of the nonempty set $A=\{nx: n\in {\cal N}\}.$ Let $s={\rm lub}(A).$ Now, there exists $z\in A$ such that $s-x<z.$ That is, there exists $m\in {\cal N}$ such that $s-x<mx.$ So that $s<(m+1)x.$ This contradicts the fact that $s$ is an upper bound of $A.$ \hfill $\Box$

The following corollary to the Archimedean principle helps in defining the integral value function in ${\cal R}.$  

\btm 
Corresponding to each real number $x$ there exists a unique integer $n$ such that $n\leq x <n+1.$ 
\etm 

\noi {\em Proof.} Let $x\in {\cal R}.$ If $x\in {\cal Z},$ then we take $n=x.$ 

If $x>0$ and $x\not\in {\cal Z},$ the Archimedean principle implies that there exists an $m\in {\cal N}$ such that $m>x.$ Let $k$ be the least natural number such that $k>x;$ so that $k-1 < x <k.$ 

If $x<0$ and $x\not\in {\cal Z},$ then $-x>0$ and $-x\not\in {\cal Z}.$ By what we have just shown, there exists a natural number $j$ such that $j-1 < -x < j.$ Then $-j < x < -(j-1)=-j+1.$ 

For uniqueness suppose that $n-1\leq x <n$ and $m-1\leq x <m$ for $n,m\in {\cal Z}.$ Assume that $n<m.$ Then $n-1<m-1$ so that we have $n-1 < m-1 \leq x < n < m.$ It implies that $m-1$ is an integer that lies between two consecutive integers $n-1$ and $n,$ and $m-1$ is not equal to either of them. This is impossible. Similarly, $m<n$ leads to a contradiction. Therefore, $n=m.$  \hfill $\Box$

Thus, the correspondence $x\mapsto n,$ where $n\leq x < n+1,$ defines a function from ${\cal R}$ to ${\cal Z}.$ We call it the {\bf integral value function}, and write it as $[x].$ That is, 
$$[x]=\h{~the~largest~integer~less~than~or~equal~to~}x\q\h{for~}x\in {\cal R}.$$ 

Once integral value function is defined, we can have a decimal representation of elements of ${\cal R}.$ 

Let $x\in {\cal R}.$ If $x=0,$ then its decimal representation is itself. 

Suppose $x>0.$  Take $x_0=[x].$ Then $0\leq y_0=x-x_0 < 1.$ Now, $0\leq 10 y_0\leq 9.$ Take $x_1=[10 y_0].$ With $y_1=y_0-x_1,$ we have $0\leq 10 y_1\leq 9.$ Take $x_2=[10y_1].$ Continuing  this process $n$ times for any $n\in\N$ gives rise to an $n$-places decimal approximation of $x$ such as 
$$x_0+\frac{1}{10}x_1+\frac{1}{10^2}x_2+\cdots + \frac{1}{10^n} x_n= x_0.x_1x_2\cdots x_n.$$ 
This process leads us to think of elements of ${\cal R}$ being placed on a straight line. We think of a positive element $x$ of ${\cal R}$ as a point on this line to the right of the point marked $0.$ To approximate it, we take the largest integer less than or equal to $x;$ this is $x_0.$ Next, we divide the line segment $x_0$ to $x_0+1$ into ten parts. If $x$ lies in the $i$th interval, then $x_1$ is this $i.$ This process continues to obtain the $n$th digit $x_n$ for each $n.$ Notice that at any stage $n,$ $~|x-x_0.x_1x_2\cdots x_n|\leq 10^{-n}.$ For each $n,$  $x$ is represented by $x_0.x_1x_2\cdots x_n$ up to $n$ decimal places. 

In case the digit $9$ repeats after a certain stage, we need to update the decimal representation. The decimal representation so obtained is in the form $a_0.a_1a_2\cdots,$ where $a_0$ is a non-negative integer and other $a_i$s are digits. If $a_i=9$ for all $i\geq 1,$ then we replace the infinite decimal with the natural number $a_0+1.$ If $a_i=9$ for all $i>m\geq 1,$ then we replace the infinite decimal with the terminating decimal $a_0.a_1\cdots a_{m-1}b,$ where $b=a_m+1.$ Leaving these cases, all other infinite decimals are kept as they are. The updated decimal, whether terminating or infinite, is referred to as the decimal representation of $x.$ 

Next, suppose $x<0.$ Then $x=-y$ for some $y>0.$  If $y_0.y_1y_2\cdots$ is the decimal representation of $y,$ then the decimal representation of $x$ is $-y_0.y_1y_2\cdots.$ 

As we know, if $x\in {\cal Q},$ then such a decimal representation recurs (possibly with all $x_i=0$ after a certain stage). If $x\not\in {\cal Q},$ then such a decimal representation is not recurring. 

Thus the decimal representation defines a function from ${\cal R}$ to $\R.$ We define  the function 
$$\phi:{\cal R}\to \R~\h{~by~}~\phi(x)=x_0.x_1x_2\cdots\in\R~\h{~for~}~x\in {\cal R}.$$ 
We will see that $\phi$ is an isomorphism. Towards this, we start with the following result. 

\btm 
Let $x,y\in {\cal R}.$ Then,  $x<y$ iff $\phi(x)<\phi(y).$ 
\etm 

\noi {\em Proof.} First, we look at the case where $0<x<y.$ We use the obvious facts  such as $\phi(n)=n$ for every $n\in {\cal N},$ and $0<1<2<\cdots <9,$ which are used in defining $<$ on $\R.$  Then $\phi(x)<\phi(y)$ follows from the construction of $\phi(z)$ for any $z\in {\cal R}$ with $0<z.$ 

Next, suppose $x<0<y.$ Then $\phi(x)$ starts with a minus sign whereas $y$ does not. So, $\phi(x)<\phi(y).$ 

Finally, if $x<y<0,$ then $0<-y < -x.$ In this case, $\phi(-y)<\phi(-x).$ As it is defined in $\R,$ $-\phi(-x)<-\phi(-y).$ However, $-\phi(-x)=\phi(x)$ and $-\phi(-y)=\phi(y).$ Therefore, $\phi(x)<\phi(y).$ This proves that if $x<y,$ then $\phi(x)<\phi(y).$ 

For the converse, notice that  if $x=y,$ then $\phi(x)=\phi(y).$ From what we have just proved, It follows that if $x>y,$ then $\phi(x)>\phi(y).$ This proves that if $x\not<y,$ then $\phi(x)\not<\phi(y).$ \hfill $\Box$

\btm 
Let $A$ be a nonempty subset of ${\cal R}$ bounded above. Then $\{\phi(x): x\in A\}$ is a nonempty subset of $\R$ bounded above, and $\phi\big({\rm lub}(A)\big) = \sup\{\phi(x): x\in A\}.$ 
\etm 

\noi {\em Proof.} Clearly, $\{\phi(x): x\in A\}\neq\empset.$ Write $\alpha={\rm lub}(A).$ For each $x\in A,~x\leq \alpha,$ and if $\beta\in \mathcal{R}$ is such that for each $x\in A,~x\leq \beta,$ then $\alpha\leq \beta.$ By the previous theorem it follows that for each $x\in A,~\phi(x)\leq \phi(\alpha),$ and if $\beta\in {\cal R}$ is such that for each $x\in A,~\phi(x)\leq \phi(\beta),$ then $\phi(\alpha)\leq \phi(\beta).$  Therefore, $\{\phi(x): x\in A\}$ is bounded above in $\R$ and $\phi(\alpha)=\sup\{\phi(x): x\in A\}.$  \hfill $\Box$

\btm\label{lubsup} 
Let $x\in {\cal R}.$ Then $x={\rm lub}\{t\in {\cal R}: t<x\}={\rm lub}\{t\in {\cal R}: t<x,~\phi(t)\in\T\}.$ 
\etm 

\noi {\em Proof.} Write $A=\{t\in {\cal R}: t<x\}.$ Now, $x$ is an upper bound of $A,$ so ${\rm lub}(A) \leq x.$ If $x\neq{\rm lub}(A),$ then ${\rm lub}(A)<x.$ So, there exists $\alpha\in {\cal R}$ such that ${\rm lub}(A)<\alpha <x.$ Then $\alpha\in A,$ and as an upper bound of $A,$ ${\rm lub}(A)<\alpha.$ This is a contradiction.  Hence $x={\rm lub}\{t\in {\cal R}: t<x\}.$ 

For the second equality, let $B=\{t\in {\cal R}: t<x,~\phi(t)\in\T\}.$ Clearly, $x$ is an upper bound of $B.$ If $x\neq {\rm lub}(B),$ then there exists $z\in B$ such that ${\rm lub}(B)\leq z <x.$ Then $\phi(z)\in\T$  and $\phi(z)<\phi(x).$ Then $\phi(x)\neq \sup\{w\in\T: w<\phi(x)\}.$ This contradicts Theorem~\ref{forall-exists-sup}(4). 
\hfill $\Box$ 

\btm 
$\phi:{\cal R}\to\R$ is a bijection. 
\etm 

\noi {\em Proof.} Since $x<y$ implies that $\phi(x)<\phi(y)$ it follows that $\phi$ is one-one. To show that $\phi$ is an onto map, we first consider all terminating decimals, and then lift the result to real numbers. 

Let $t=\pm a_1\cdots a_k.b_1b_2\cdots b_m\in\T.$ Write $n=\pm a_1\cdots a_kb_1b_2\cdots b_m$ and $x=10^{-m}\times n.$ Consequently, $n\in {\cal Z}$ and $x\in {\cal R}.$ We see that $\phi(x)=t.$ 

Next, let $t\in \R\setminus\T.$ Then $t=\sup\{s\in \T:  s<t\}.$ By what we have just proved, each $s\in\T$ is $\phi(u)$ for some $u\in {\cal R}.$ Thus  $t=\sup\{\phi(u): u\in {\cal R},~\phi(u)\in\T,~\phi(u)<t\}.$ By Theorem~\ref{lubsup}, $t=\phi\big({\rm lub}\{u\in {\cal R}: \phi(u)\in\T, \phi(u)<t\} \big).$ Then, with $x={\rm lub}\{u\in {\cal R}: \phi(u)\in\T, \phi(u)<t\} ,$ we see that $t=\phi(x).$ 

Therefore, $\phi$ is an onto map. \hfill $\Box$

\btm 
Let $x,y\in {\cal R}.$ Then $\phi(x+y)=\phi(x)+\phi(y).$ 
\etm 

\noi {\em Proof.} We first show the result for those elements in ${\cal R}$ whose decimal representations are terminating decimals. For any integer $n$ in ${\cal R},$ $\phi(n)=n.$  Thus $\phi(x+y)=\phi(x)+\phi(y)$ holds for integers $x$ and $y.$ Then by multiplying a suitable negative power of $10,$ the same equality is proved for terminating decimals. The details follow. 

Let $s,t\in {\cal R}$ be such that $\phi(s),\phi(t)\in\T.$ If the number of digits after the decimal points in $\phi(s)$ and $\phi(t)$ are unequal, then adjoin the necessary number of zeros to one of them so that the number of digits after the decimal points in both are same. So, let $\phi(s)=\pm s_1\cdots s_j.a_1a_2\cdots a_m$ and $\phi(t)=\pm t_1\cdots t_k.b_1b_2\cdots b_m.$ Then $s=10^{-m}\times (\pm s_1\cdots s_ja_1a_2\cdots a_m),$ $t= 10^{-m}\times (\pm t_1\cdots t_kb_1b_2\cdots b_m)$ and $s+t=10^{-m}\times (\pm c_1\cdots c_id_1d_2\cdots d_m)$ for suitable digits $c$s and $d$s. It follows that $$\phi(s+t)=\pm c_1\cdots c_i.d_1d_2\cdots d_m =\phi(s)+\phi(t).$$ 
Here, the signs $\pm$ are taken suitably. 

Next, we show the result for those elements of ${\cal R}$ whose decimal representations are  non-terminating decimals. Towards this, let $r,s,t,x,y\in {\cal R}.$ We observe that if $r<x$ and $s<y,$ then $r+s<x+y.$ Further, if $t<x+y,$ then write $x+y-t=\epsilon$ with some $\epsilon>0.$ Take $r=x-\epsilon/2$ and $s=y-\epsilon/2.$ Then $r+s=x+y-\epsilon =t.$ That is, if $t<x+y,$ then there exist $r,s$ such that $r<x,~s<y$ and $t=r+s.$ We use these observations in the following calculation. 

Let $x,y\in {\cal R}$ be such that $\phi(x),\phi(y)\in\R\setminus\T.$ Then 
\beqarray 
\phi(x+y) & = & \phi\big({\rm lub}\{t\in {\cal R}: t<x+y,~\phi(t)\in\T\} \big) \\  
& = & \sup\{\phi(t): t\in {\cal R},~t< x+y,~ \phi(t)\in\T\} \\ 
& = & \sup\{\phi(r+s): r,s\in {\cal R},~r+s<x+y,~\phi(r+s)\in\T \} \\ 
& = & \sup\{\phi(r)+\phi(s): r,s\in {\cal R},~r+s<x+y,~\phi(r+s)\in\T \} \\ 
& = & \sup\{\phi(r)+\phi(s): r,s\in {\cal R},~r<x,~s<y,~\phi(r),\phi(s)\in\T \} \\ 
& & {\rm (It~follows~from~the~above~observations~and~Theorem~\ref{forall-exists-sup}.)}  \\ 
& = & \sup\{\phi(r)+\phi(s): r,s\in {\cal R},~\phi(r)<\phi(x),~\phi(s)<\phi(y),~\phi(r),\phi(s)\in\T \} \\ 
& = & \sup\{u+v: u<\phi(x),~v<\phi(y),~u,v\in\T \}  =  \phi(x)+\phi(y). \hspace{2cm} \Box 
\eeqarray

\btm 
Let $x,y\in {\cal R}.$ Then $\phi(x\times y)=\phi(x)\times\phi(y).$ 
\etm 

\noi {\em Proof.} Clearly, if one of $x$ or $y$ is equal to $0,$ the equality holds. For $x,y\in {\cal N},$  $\phi(x)=x,~\phi(y)=y$ and $\phi(x+y)=x+y;$  thus the equality holds. Next, suppose $x,y\in {\cal R}$ are such that $\phi(x),\phi(y)$ are positive terminating decimals. Then $x>0$ and $y>0.$ As in the case of addition, we see that $x=10^{-m}k_1$ and $y=10^{-n}k_2$ for $n_1,n_2\in N.$ Then $x\times y=10^{-m-n}(k_1\times k_2)$ so that the equality holds. 

Next, suppose that $x,y\in {\cal R}$ are such that $\phi(x)$ and $\phi(y)$ are positive non-terminating decimals. Before proceeding towards the equality, we observe the following: 

Let $r,s,t,x,y\in {\cal R},$ where all of these are greater than $0.$ If $r<x$ and $s<y,$ then $r\times s<x\times y.$ Further, if $t<x\times y,$ write $t/(x\times y)=\epsilon.$ Then $\epsilon<1.$ Take $r=x\times \sqrt{\epsilon}$ and $s=y\times \sqrt{\epsilon}.$ Then $r\times s=x\times y\times \epsilon =t.$ That is, if $t<x\times y,$ then there exist $r,s$ such that $r<x,~s<y$ and $t=r\times s.$ 

Using these observations a calculation similar to that in the case of addition can now be given, where we replace all occurrences of $+$ with $\times$ to obtain the required equality. 

When one or both of $x,y$ is (are) less than $0,$ we use $(-x)\times y=-(x\times y),$ etc. for proving the equality. 
\hfill $\Box$ 

To summarize, we have proved the following for any $x,y\in {\cal R}$ and any nonempty subset $A$ of ${\cal R}$: 
\begin{enumerate}   
\item If $x<y,$ then $\phi(x)<\phi(y).$ 
\item $\phi(x+y)=\phi(x)+\phi(y).$ 
\item $\phi(x\times y)=\phi(x)\times \phi(y).$ 
\item If $A$ is bounded above, then $\{\phi(x): x\in A\}$ is a nonempty subset of $\R$ bounded above; and  $\phi({\rm lub}(A))=\sup\{\phi(x): x\in A\}.$ 
\end{enumerate}  

That is, the following result has been proved. 

\btm 
The function $\phi: {\cal R}\to\R$ given by $\phi(x)=$ the decimal representation of $x,$ is an isomorphism. 
\etm 

\section{Conclusions} 
Now onwards, we do not distinguish between ${\cal N}$ and $\N,$ between ${\cal Z}$ and $\Z,$ between ${\cal Q}$ and $\Q,$ and between ${\cal R}$ and $\R.$ We use the latter symbols instead of the former. That is, we have a unique complete ordered field, namely, $\R$ and it contains $\N,~\Z$ and $\Q.$ Also, each number in $\R,$ called a real number, is an infinite decimal. The rational numbers are those which can be written in the form $p/q$ for $p\in\Z$ and $q\in\N,$ where $p,q$ have no common factors. Again, a rational number is an infinite decimal where after some finite number of digits, some finite sequence of digits is repeated infinitely often. As usual, a terminating decimal is identified with an infinite decimal with trailing $0$s and no infinite decimal ends with a sequence of the digit $9.$ We have seen that the supremum or least upper bound of a nonempty set of real numbers that is bounded above, is unique. In accordance with this, we do not distinguish between $\sup$ and ${\rm lub}$ of such a set.  The real numbers which are not rational numbers, called the irrational numbers, are non-recurring infinite decimals. For instance, $\sqrt{2}$ is an irrational number. 

We see that the usual intuition of defining a real number as an infinite decimal is viable and workable. The theory of infinite decimals not only provides a model for a complete ordered field, but it defines real numbers in a unique manner, up to an isomorphism. The modern definition of the real number system as a complete ordered field is consistent relative to the system of natural numbers, and is categorical. Therefore, all constructions of real numbers as pointed out in \cite{weis} construct the same object $\R.$

\end{document}